\catcode`\@=11

\font\tenmsb=msbm10
\font\sevenmsb=msbm7 \font\fivemsb=msbm5  \newfam\msbfam
\textfont\msbfam=\tenmsb\scriptfont\msbfam=\sevenmsb
\scriptscriptfont\msbfam=\fivemsb

\def\hexnumber@#1{\ifnum#1<10 \number#1\else \ifnum#1=10 A\else\ifnum#1=11
 B\else\ifnum#1=12 C\else \ifnum#1=13 D\else\ifnum#1=14 E\else\ifnum#1=15
 F\fi\fi\fi\fi\fi\fi\fi}
 \def\msb@{\hexnumber@\msbfam}

\mathchardef\hbar="0\msb@7E

\def\Bbb{\ifmmode\let\next\Bbb@\else
\def\next{\errmessage{Use \string\Bbb\space only in math mode}}\fi\next}
\def\Bbb@#1{{\Bbb@@{#1}}} \def\Bbb@@#1{\fam\msbfam#1}

\catcode`\@=\active
\def\CL{\hbox{{$\cal L$}}}
\def\cg{\hbox{{\sl g}}} % used for Lie algebra 'gothic g'

\def\C{{\Bbb C}}
\def\Z{{\Bbb Z}}

\def\extd{{\rm d}}
\def\del{{\partial}}
\def\eps{{\epsilon}}
\def\tens{\mathop{\otimes}}

\def\und#1{{\underline{#1}}}
\def\equad{\kern -1.7em}

\def\eqn#1#2{\begin{equation}#2\label{#1}\end{equation}}

\def\cmath#1{\[\begin{array}{c} #1 \end{array}\]}

\def\ceqn#1#2{\begin{equation}\label{#1}
\begin{array}{c}#2\end{array}\end{equation}}

\def\note#1{}

\documentstyle[11pt,epsf]{article}
\begin{document}

\newtheorem{lemma}{Lemma}[section]
\newtheorem{propos}[lemma]{Proposition}
\newtheorem{example}[lemma]{Example}
\newtheorem{theorem}[lemma]{Theorem}
\newtheorem{cor}[lemma]{Corollary}
\newtheorem{corol}[lemma]{Corollary}
\newtheorem{defin}[lemma]{Definition}
\newtheorem{remark}[lemma]{Remark}
\newtheorem{conjec}[lemma]{Conjecture}

{\ }\qquad\qquad \hskip 4.3in  DAMTP/97-63
 \vspace{-.2in}

\begin{center} {\LARGE ANYONIC LIE ALGEBRAS}
\\  {\ }
{\ }\\ S. Majid\footnote{Royal Society University Research Fellow
and Fellow of Pembroke College, Cambridge, England.\\
This paper is in final form and no version of it will be published
elsewhere.} \\ {\ }\\
Department of Applied Mathematics \& Theoretical Physics\\
University of Cambridge, Cambridge CB3 9EW, UK\\
www.damtp.cam.ac.uk/user/majid\\
\end{center}

\vspace{10pt}
\begin{quote}\baselineskip 13pt
\noindent{\bf Abstract}
 We introduce anyonic Lie algebras in terms of structure constants. We
 provide the simplest examples and formulate some open problems.
\end{quote}

\baselineskip 22pt

\section{Introduction}

There have been a number of attempts to generalise supersymmetry in
physics to anyonic or fractional statistics. The approach described
here originates\cite{Ma:any} in the modern theory of quantum groups
and braided groups. Although not yet reaching the stage of
Lagrangians and field theory, it is mathematically well founded and
could perhaps be incorporated into future theories. After a brief
introduction to this approach, we describe a natural notion of
`anyonic Lie algebra' arising as a special case of \cite{Ma:lie},
but explicitly formulated now through structure constants obeying
certain conditions. We also formulate the classification and other
problems, and provide some simple examples with $\Z_{/3}$-grading.

\section{Anyspace}

The easiest anyonic object to understand is 1-dimensional anyspace.
This is the algebra $\C[\theta]/\theta^n$ with one coordinate
$\theta$ obeying $\theta^n=0$ and of degree $|\theta|=1$. $n=2$ is
a usual Grassmann variable. One may add anyonic variables by
\eqn{anyadd}{ \theta''=\theta+\theta',\quad \theta'\theta
=e^{2\pi\imath\over n}\theta\theta'}
where $\theta'$ is another copy if the anyonic variable with
statistics as shown. One may check that $\theta''^n=0$ using
properties of $q$-binomial coefficients. This additivity provides
the key properties analogous to the real line. It has been
introduced in \cite{Ma:any}.

To be precise one should use the notion of `braided
groups'\cite{Ma:bra} or Hopf algebras with braid statistics.
Firstly, a braided group is an object in a braided category -- a
category for which the `exchange' or braiding between any two
objects is coherently specified. In our case we use  the  category
of $n$-anyonic vector spaces where objects are $\Z_{/n}$-graded
spaces and the `braided transposition' is
\eqn{anybraid}{\Psi(x\tens
y)=e^{2\pi\imath |x||y|\over n}y\tens x} on elements $x,y$ of
homogeneous degree $|\ |$. This category has been introduced in
\cite{Ma:any} as the category of representations of a certain
quantum group $\Z_{/n}'$. Morphisms preserve degree, although one
also (as in SUSY) considers linear maps which are not
degree-preserving. Secondly, the group law is expressed as a
morphism
\eqn{anydelta}{ \Delta:\C[\theta]/\theta^n\to
\C[\theta]/\theta^n\ \und\tens\ \C[\theta]/\theta^n,\quad
\Delta\theta=\theta\tens 1+1\tens\theta.}
If we write $\theta=\theta\tens 1$ and $\theta'=1\tens \theta$
(i.e. the generator of the first copy of the anyonic line is called
$\theta$ and the generator of the second copy is called $\theta'$)
then we meet the previous notation. We underline $\und\tens$
because the two copies do not commute -- they have the {\em braided
tensor product algebra} $(x\tens y)(w\tens z)=x\Psi(y\tens w)z$. In
the formal specification of a braided group we must provide also a
counit $\eps(\theta)=0$ and an antipode $S(\theta)=-\theta$.

In general, a braided group means an algebra $B$ with braiding
$\Psi:B\tens B\to B\tens B$, a coproduct $\Delta:B\to B\tens B$,
counit $\eps:B\to\C$ and antipode $S:B\to B$. The easiest way to
write the axioms is in a diagrammatic notation\cite{Ma:bra}, see
Figure~1, in which we write all maps pointing generally downwards, with
$\cdot=\epsfbox{prodfrag.eps}$, $\Delta=\epsfbox{deltafrag.eps}$,
$\Psi=\epsfbox{braid.eps}$, $\Psi^{-1}=\epsfbox{braidinv.eps}$.
Other morphisms are nodes with the appropriate number of inputs and
outputs. In this notation we `wire' the outputs of maps into the
inputs of other maps to construct our algebraic operation.
Information flows along these wires much as in a computer, except
that under and over crossings are nontrivial operators
$\Psi,\Psi^{-1}$. This is a new kind if `braided mathematics'. Note
that physicists use such `wiring' notation in the form of Feynman
diagrams. In anyonic or braided field theory one would have
such diagrams but with under or over crossings being nontrivial
operators.
\begin{figure}
\[ \epsfbox{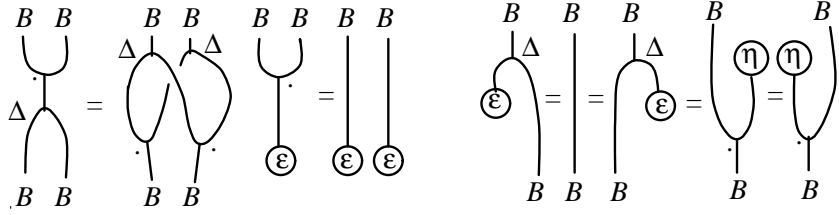}\]
\caption{Homomorphism and counity axioms for a braided group}
\end{figure}
The diagrams translate into equations by reading the operations
from the top down. Thus an anyonic braided group means a
$\Z_{/n}$-graded algebra $B$ and coalgebra defined by $\eps:B\to
\C$ and $\Delta(x)=x_A\tens x^A$ say (summation understood over terms
labeled by $A$ of homogeneous degree), coassociative and counital
in the sense
\eqn{coass}{x_{AB}\tens x_A{}^B\tens x^A=x_A\tens x^A{}_B
\tens x^{AB},\quad \eps(x_A)x^A=x=x_A\eps(x^A)}  and obeying
\eqn{anybg}{ (xy)_A\tens (xy)^A=x_A x_B\tens x^A x^B e^{{2\pi\imath\over n}
|x^A||x_B|}}
for all $x,y\in B$. The axioms for the antipode are as for usual quantum
groups\cite{Ma:book} and we omit writing them explicitly.

One has many further
structures, such as integration, differentiation
\eqn{anyint}{ \int\theta^m=\cases{1&m=n-1\cr0&else},\quad\del f(\theta)
={f(\theta)-f(e^{2\pi\imath\over n}\theta)\over(1-e^{2\pi\imath\over n})\theta}}
as well as exponentials, Gaussians\cite{BauFla:pat} and
$\delta$-functions\cite{MaPla:ran}. In the paper\cite{MaPla:ran}
one finds also random walks and Brownian motion on anyspace.

Higher-dimensional anyonic planes can be defined by tensor product.
However, if $B,C$ are braided groups their natural braided tensor
product algebra and coalgebra $B\und\tens C$ does not form a
braided group in the original category. In the diagrammatic
notation, one gets `tangled up' if one tries this. Instead, one has
to `glue' the categories also\cite{MaMar:glu}. This is known for
Hecke-type braid statistics where the braiding has two eigenvalues.
 In our setting it means we can `glue'
anyonic variables of degree either 1 (as for $\theta$ above) or
$n/2-1$ when $n$ is even. We call the latter type of coordinate
$\eta$. So a natural higher-dimensional anyspace has an
$s$-dimensional glued anyonic part $\{\theta^i\}$ and an
$r$-dimensional glued `fermi-anyonic' part $\{\eta^i\}$. The
required `glued' braiding for additivity is, however, no longer
just a phase
-- it involves linear combinations (actually it is the non-standard
$sl_{s|r}$ R-matrix as explained in \cite{MaMar:glu}.) Infinite
anyonic planes can be found among primary fields in 2-d quantum
gravity as explained in \cite{Ma:inf}. In short, the `phase-factor'
anyonic form (\ref{anybraid}) is not closed under tensor product
and leads one naturally into systems with linear combinations in
the statistics!

Although not built by tensor product, there are still plenty of
higher-dimensional anyonic braided groups. For example, one has
anyonic matrices\cite{MaPla:any} as a generalization of quantum
matrices \cite{FRT:lie} and supermatrices. Here there are $N^2$
generators $\{t^i{}_j\}$ of degree $|t^i{}_j|=f(i)-f(j)$ where $f$
is a degree in $\Z_{/n}$ associated with the row or column and
\eqn{anymatcop}{e^{{2\pi\imath\over n}(f(i)f(k)+f(j)f(b))}
{\und R}^i{}_a{}^k{}_b t^a{}_j t^b{}_l
=e^{{2\pi\imath\over n}(f(j)f(l)+f(i)f(b))}t^k{}_b t^i{}_a
{\und R}^a{}_j{}^b{}_l,
\quad \Delta t^i{}_j=t^i{}_a\tens t^a{}_k,\quad \eps t^i{}_j
=\delta^i{}_j.}
We require that $\und R$ obeys certain anyonic-Yang-Baxter
equations\cite{MaPla:any} (with the result that the anyonic
matrices are any-coquasitriangular). One method to obtain $\und R$
is to start with certain solutions $R$ of the usual braid or
Yang-Baxter equations and `transmute' them. The quotients of such
anyonic matrices yield further anyonic braided groups, see
\cite{MaPla:any}\cite[Appendix]{Ma:varen}. However, the method
probably does not exhaust all anyonic braided groups.

\bigskip
{\em Problem 1} Classify all low-dimensional anyonic braided
groups, i.e find all algebras, $\Delta,\eps$ and optionally
antipodes $S$ for a given $\Z_{/n}$-graded vector space $B$.
\bigskip

For $n=1$ (usual quantum groups) the low-dimensional classification
has been done by Radford via computer. There are also some general
techniques \cite{LarRad:fin}. For $n=2$ (super quantum groups) one has
many examples and probably similar techniques. $n=3$ would be
the simplest truly braided case, with $\Z_{/3}$-anyonic braiding
and is wide open. Examples are in \cite{Ma:any}\cite{MaPla:any}.

\section{Anyonic Lie algebras}

The axioms of a braided Lie-algebra are shown in Figure~2 in a
diagrammatic form\cite{Ma:lie}\cite{Ma:sol}. We need an object
$\CL$, a morphism $[\ ,\ ]$, and (unusually) morphisms
$\Delta:\CL\to\CL\tens\CL$ and $\eps:\CL\to
\C$ forming a coalgebra. While there are one or two recent
proposals for a quantum or braided Lie algebra, this is the only
one which (a) includes quantum groups and (b) has features of Lie
algebras such as tensor product representations, enveloping
(braided) bialgebra, etc.  We now specialize this theory to the
anyonic setting.
\begin{figure}
\[ \epsfbox{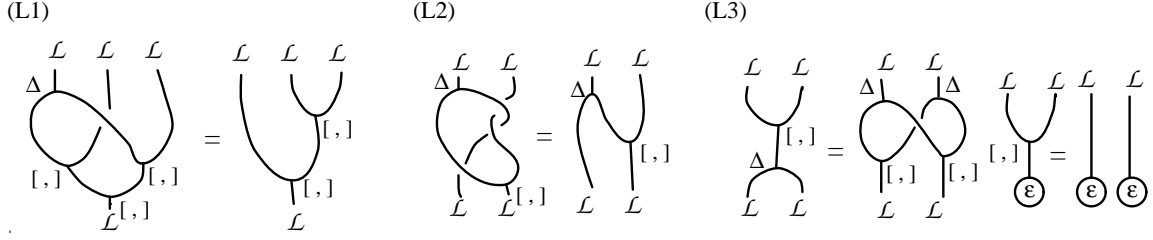}\]
\caption{Axioms of a braided-Lie algebra}
\end{figure}

Thus, an anyonic Lie algebra means a $\Z_{/n}$-graded vector space
$\CL$, say, and degree-preserving maps $\Delta,\eps$ obeying
(\ref{coass}) and

\eqn{liea}{\eps([x,y])=\eps(x)\eps(y),\quad [x,y]_A\tens [x,y]^A=[x_A,y_B]
\tens[x^A,y^B]e^{2\pi\imath |x^A||y_B|\over n}}
\eqn{lieb}{x_A\tens
[x^A,y]=x^A\tens [x_A,y]
e^{{2\pi\imath\over n}|x^A|(2|y|+|x_A|)}}
\eqn{liec}{ [x,[y,z]]=[[x_A,y],[x^A,z]]
e^{2\pi\imath |y| |x^A|\over n}.} This is obtained by reading off
the diagrams with the anyonic braiding (\ref{anybraid}). We write
$\Delta x=x_A\tens x^A$ as a notation (summation over labels $A$
understood). The anyonic enveloping algebra $U(\CL)$ comes out from
the diagrammatic definition in \cite{Ma:lie} as generated by
products of $\CL$ with the relations
\eqn{lied}{ xy=[x_A,y]x^A e^{2\pi\imath |x^A||y|\over n}}
and $\Delta,\eps$ extended as an anyonic braided group by
(\ref{anybg}).

If we fix a basis $\{x^\mu\}$ of degree $|x^\mu|=p(\mu)\in\Z_{/n}$,
then these maps are equivalent to structure constants $\eps,d,c$
defined by
\eqn{anystruc}{ \eps(x^\mu)=\eps^\mu,\quad \Delta x^\mu
=d^\mu{}_{\nu\rho}x^\nu\tens x^\rho,\quad
[x^\mu,x^\nu]=c^{\mu\nu}{}_\rho x^\rho} obeying
\eqn{liecompa}{ \eps^\mu p(\mu)=0,
\quad d^\mu{}_{\nu\rho}(p(\mu)-p(\nu)-p(\rho))=0,
\quad c^{\mu\nu}{}_\rho(p(\mu)+p(\nu)-p(\rho))=0}
\eqn{liecompb}{ d^\mu{}_{\alpha\lambda}d^\alpha{}_{\nu\rho}
=d^\mu{}_{\nu\alpha}d^\alpha{}_{\rho\lambda},
\quad d^\mu{}_{\alpha\nu}\eps^\alpha=\delta^\mu{}_\nu
=d^\mu{}_{\nu\alpha}\eps^\alpha,\quad c^{\mu\nu}{}_\alpha \eps^\alpha
=\eps^\mu \eps^\nu}
\eqn{liecompc}{c^{\mu\nu}{}_\alpha d^\alpha{}_{\rho\lambda}
=e^{{2\pi\imath\over n} p(\beta)p(\gamma)}d^\mu{}_{\alpha\beta}
d^\nu{}_{\gamma\delta}c^{\beta\delta}{}_\lambda
c^{\alpha\gamma}{}_\rho}
\eqn{liecompd}{ d^\mu{}_{\lambda\alpha}
c^{\alpha\nu}{}_\rho=e^{{2\pi\imath\over n} p(\lambda)(2p(\nu)
+p(\alpha))}d^\mu{}_{\alpha\lambda}c^{\alpha\nu}{}_\rho}
\eqn{liecompe}{ c^{\mu\nu}{}_\alpha c^{\rho\alpha}{}_\lambda
=e^{{2\pi\imath \over n}p(\beta)p(\mu)}d^\rho{}_{\alpha\beta}
c^{\alpha \mu}{}_\gamma c^{\beta \nu}{}_\delta c^{\gamma \delta}{}_\lambda.}
The associated anyonic braided group $U(\CL)$ is generated by $1,x^\mu$ with
the quadratic relations
\eqn{liecompf}{ x^\mu x^\nu=e^{{2\pi\imath\over n}p(\beta)p(\nu)}
d^\mu{}_{\alpha\beta}c^{\alpha \nu}{}_\gamma x^\gamma x^\beta}

A particular ansatz is of the form $\CL=\C\oplus \cg$, where
$\C$ is spanned by $x^0$ with $p(0)=0$ and $\cg$ is spanned by the
remaining $x^i$, $i\ge 1$, and
\ceqn{ansatz}{  \eps^0=1,\quad \eps^i=0,\quad d^0{}_{00}=1,
\quad d^0{}_{ij}=d^0{}_{i0}=d^0{}_{0i}=d^i{}_{00}=d^i{}_{jk}=0,
\quad d^i{}_{0j}=d^i{}_{j0}=\delta^i{}_j  \\
 c^{00}{}_0=1,\quad c^{00}{}_i=c^{0i}{}_0=c^{i0}{}_0=c^{i0}{}_j
 =c^{ij}{}_0=0,\quad c^{0i}{}_j=\delta^i{}_j.}
In this case (\ref{liecompa})-(\ref{liecompe}) for the remaining
variables $p(i), c^{ij}{}_k$ reduce to
\eqn{liesuper}{ 1=e^{{2\pi\imath\over n}2p(i)p(j)},\quad c^{ij}{}_a
c^{ka}{}_l=c^{ki}{}_ac^{aj}{}_l+e^{{2\pi\imath
\over n}p(k)p(i)}c^{kj}{}_a c^{ia}{}_l}
from (\ref{liecompd}) and (\ref{liecompe}) respectively. This means
effectively that $n=1,2$ and the ansatz recovers the Jacobi axiom
for a usual or super-Lie algebra (other $n$ are equivalent after
redefining $p(i)$). The associated anyonic braided group $U(\CL)$
is generated by $1,x^0,x^i$ with the quadratic relations
\eqn{liesuperenv}{ x^0x^i=x^ix^0,\quad x^ix^j=e^{{2\pi\imath\over n}
p(i)p(j)}x^jx^i + c^{ij}{}_a x^a x^0}
which is the homogenised usual or super enveloping algebra. Setting
the central element $x^0=1$ recovers the usual $U(\cg)$. Actually,
our objects even for $n=1,2$ are slightly more general than usual
or super-Lie algebras because we do not impose any equation like
\eqn{lieanti}{{}[x^i,x^j]=-e^{{2\pi\imath\over n}
p(i)p(j)}[x^j,x^i].} It turns our that one does not need such
antisymmetry for the most important parts of Lie theory or
super-Lie theory, although one can impose it additionally at least
in the context of the ansatz (\ref{ansatz}). I do not know a
general diagrammatic or axiomatic way to introduce such an
antisymmetry condition, however. A related problem: the braided
enveloping algebra $U(\CL)$ does not usually have an antipode. For
this, one needs to quotient it (as in the above ansatz where we set
$x^0=1$). In all known examples the quotienting procedure is also
known, but `by hand'.

\section{Examples of Matrix Type}

We now give a family of anyonic Lie algebras $\CL_{N,f}$ going
genuinely beyond the super and usual ones. They are a
specialization of the general R-matrix construction in
\cite{Ma:lie} to  the case
\eqn{R}{ R^m{}_n{}^r{}_l=\delta^m{}_n\delta^r{}_l e^{{2\pi\imath
\over n}f(m)f(r)}}
where $m,n,r,l=1,\cdots,N$ and $f(m)\in \Z_{/n}$ is an arbitrary
grading function. We write the $N^2$-dimensional vector space
$\CL_{N,f}$ in a matrix (`twistor') form where $x^\mu=x^m{}_{\dot m
}$ and $\mu$ corresponds to the multiindex $(m,\dot m)$. Computing
from \cite[Prop. 5.2]{Ma:lie}, the induced braiding $\Psi$ comes
out to be the anyonic one with $p(\mu)=f(m)-f(\dot m)$, and
\eqn{matanycomp}{ d^\mu{}_{\nu\rho}=\delta^m{}_n\delta^{\dot n}{}_{r}
\delta^{\dot r}{}_{\dot m},\quad \eps^\mu=\delta^m{}_{\dot m},\quad
c^{\mu\nu}{}_\rho=\delta^m{}_{\dot m}\delta^n{}_r\delta^{\dot
r}{}_{\dot n} e^{-{2\pi\imath\over n}2 f(m)p(\nu)}} solves the
equations (\ref{liecompa})--(\ref{liecompe}). Equivalently,
\eqn{matanylie}{\Delta x^m{}_{\dot m}=x^{m}{}_a\tens x^a{}_{\dot m},
\quad \eps(x^m{}_{\dot m})=\delta^m{}_{\dot m},\quad [x^\mu,x^\nu]
=\eps(x^\mu)x^\nu e^{-{2\pi\imath\over n}2f(m)p(\nu)}.}
The anyonic enveloping algebra $U(\CL_{N,f})$ has the relations
\eqn{matanyenv}{ x^\mu x^\nu=e^{-{2\pi\imath\over n}(f(m)
+f(\dot m))p(\nu)}x^\nu x^\mu} and
the coproduct extended as an anyonic braided
group.

The 1-dimensional case goes naturally with the bosonic choice
$n=1$ and has the  structure
\[ [x,x]=x,\quad \Delta x=x\tens x,\quad \eps(x)=1,\quad |x|=0,
\quad U(\CL_{1,f})=\C[x].\]
It is the structure of the central element $x^0$ in the ansatz
(\ref{ansatz}) above.

The 4-dimensional case $\CL_{2,f}$ has generators
$x^\mu=\pmatrix{a&b\cr c&d}$ say. With $n=2$ and $f(1)=0, f(2)=1$
as elements of $\Z_{/2}$, we have $a,d$ bosonic and $b,c$
fermionic, and a super-anyonic Lie algebra
\cmath{{} \eps(a)=\eps(d)=1,\quad \eps(b)=\eps(c)=0,\quad
 \Delta a=a\tens a+b\tens c\\
 \Delta b=b\tens d+a\tens b,\quad \Delta c=c\tens a+d\tens c,\quad
 \Delta d=d\tens d+c\tens b\\
{}[a,x]=[d,x]=x,\quad [b,x]=[c,x]=0} for all $x=a,b,c,d$. Note that
even although the braiding is the super or $\Z_{/2}$-graded one,
this $\CL_{2,f}$ is not a usual super-Lie algebra because the
coproduct does not have the simple form as in the ansatz in the
last section; it is a new kind of super structure. On the other
hand, its enveloping super-bialgebra $U(L_{2,f})$ has relations
\[ ax=xa,\quad dx=xd,\quad b^2=c^2=0,\quad bc= -cb \]
for all $x=a,b,c,d$, i.e. $U(L_{2,f})=M_{1|1}$, the usual algebra
of $2\times 2$ super-matrices.

A different choice for $f$ is with $n=3$ and $f(1)=0, f(2)=1$ as
elements of $\Z_{/3}$, which gives $\CL_{2,f}$ with a similar
$2\times 2$ matrix of generators but with degrees
\[ |a|=|d|=0,\quad |b|=-1,\quad |c|=1\]
in $\Z_{/3}$. This gives us $\CL_{2,f}$ as an anyonic Lie algebra
with $\Z_{/3}$-grading. It has the structure
\cmath{{} \eps(a)=\eps(d)=1,\quad \eps(b)=\eps(c)=0,\quad
\Delta a=a\tens a+b\tens c\\
 \Delta b=b\tens d+a\tens b,\quad \Delta c=c\tens a+d\tens c,\quad
 \Delta d=d\tens d+c\tens b\\
{}[a,x]=x,\quad [d,x]=e^{{2\pi\imath\over 3}|x|}x,\quad
[b,x]=[c,x]=0} for all $x=a,b,c,d$. The enveloping
$\Z_{/3}$-anyonic braided group has relations
\[ ax=xa,\quad cb=e^{2\pi\imath\over 3}bc,\quad b^2=c^2=bd=db=dc=cd=0\]
and coproduct extended to the algebra via (\ref{anybg}). One may
check that $ad-cb$ is bosonic, group-like and central and hence may
be set equal to 1. Then $b=c=0$ and $a=d^{-1}$. Hence the quotient
of this $U(\CL_{2,f})$ by the determinant relation is the Hopf
algebra $\C[a,a^{-1}]$.

This second $\CL_{2,f}$ is probably the simplest truly braided
anyonic Lie algebra. We can similarly choose various $n$ and $f$
for larger matrices, such as $n=N$ and $f(i)=i-1$ in $\Z_{/n}$ (as
for the first $\CL_{2,f}$ above). The $\CL_{3,f}$ of this type has
a matrix of generators
\[ x^\mu=\pmatrix{a& b_{-}& b_+\cr c_{+}&d_{+}&e_-\cr c_-&e_{+}&d_{-}},
\quad |x^\mu|=\pmatrix{0&-1&1\cr1&0&-1\cr-1&1&0}\] with the matrix
coalgebra, and the $\Z_{/3}$-anyonic Lie bracket
\[ [a,x]=x,\quad [d_\pm,x]=e^{\pm{2\pi\imath |x|\over 3}}x,
 \quad [b_\pm,x]=[c_\pm,x]=[e_\pm,x]=0.\] The
enveloping $\Z_{/3}$-anyonic bialgebra $U(\CL_{3,f})$ has the
following form. $a$ is central, $b_\pm^2=b_\pm b_\mp=0$ and
similarly for $c_\pm$, but products of $b,c$ quasi-commute
according to
\[ b_\pm c_\pm=e^{2\pi\imath\over 3}c_\pm b_\pm,\quad
b_\pm c_\mp=e^{-{2\pi\imath\over 3}}c_\mp b_\pm.\] By contrast, the
$d_\pm$ form a commutative polynomial subalgebra, as do the
$e_\pm$, but the product of $d$ with $e$ variables is zero.  The
product of $b$ or $c$ variables with $d$ or $e$ variables is also
zero.

\bigskip
{\em Problem 2} Classify all low-dimensional anyonic Lie algebras
by finding all possible $\Delta,[\ ,\ ],\eps$ for a given grading
$p(\mu)$ for each basis element of $\CL$.

\bigskip
{\em Problem 3} Extend the axioms of a braided-Lie algebra to provide
in general for a quotient of $U(\CL)$ forming a braided group with
antipode.
\bigskip

Actually, one has a matrix braided Lie algebra over all regular
points of the `Yang-Baxter variety' as explained in \cite{Ma:lie}
with (\ref{R}) being special. Deformations of the $n=1,2$ cases are
controlled by the classical Yang-Baxter equations and lead to
quantum and super-quantum groups. Also, by letting $|x|\in G$, an
Abelian group and taking braiding $\Psi(x\tens y)
=y\tens x \beta(|x|,|y|)$ where $\beta$ is a bicharacter, one has in a
basis $\{x^\mu\}$ the similar axioms with $\beta(p(\mu),p(\nu))$ in
place of $e^{{2\pi\imath\over n}p(\mu)p(\nu)}$ in
(\ref{liecompa})--(\ref{liecompf}). This is a general class of
braided-Lie algebras which is not exactly anyonic (there need not
be a root of unity in the picture) but has a similar `phase factor'
form. The ansatz of the form (\ref{ansatz}) with $p(x^0)$ the group
identity and $\beta$ skew in the sense $\beta(g,h)
=\beta^{-1}(h,g)$ recovers the axioms of a colour Lie algebra. More
generally, we are not limited to skew bicharacters, so we
generalize the theory of colour-Lie algebras immediately (anyonic
ones being an example of this generalisation). One may similarly
use a bicharacter in (\ref{R}) for matrix braided Lie algebras
$\CL_{N,\beta}$.

More general braided Lie algebras $\cg_q$ are associated to the
quantum groups\cite{Dri} $U_q(\cg)$ as the appropriate
infinitesimal generators. For the $ABCD$ series we use the R-matrix
construction from \cite[Prop.~5.2]{Ma:lie} where
$\cg_q=\{x^i{}_j\}$ and $U(\cg_q)$ comes out as the braided matrix
algebra
\eqn{BR}{R_{21}x_1Rx_2=x_2R_{21}x_1R}
(in a compact notation) for the appropriate R-matrix (this is
\cite[eqn. (20)]{Ma:exa} for the braided matrix algebra $B(R)$,
written compactly). For the $A$ series one has $\cg_q$ as a
deformation of the $\CL=u(1)\oplus\cg$ ansatz (\ref{ansatz}).
Similarly to that, the corresponding $U(\cg_q)$ has to be
quotiented by setting a Casimir=1 to obtain (some version of) the
Drinfeld-Jimbo quantum enveloping algebras $U_q(\cg)$ by this
method.  To explain this for the general case, we make a
`quantum-geometry transformation'\cite{Ma:mex} and indeed view the
$\{x^i{}_j\}$ as  braided matrix `coordinates' of $N\times N$
braided matrices $B(R)$\cite{Ma:exa}. Since these were introduced
(by the author) as analogues of the quantum matrices $A(R)$ in
\cite{FRT:lie}, they similarly  have quotients by (braided)
$q$-determinant=1 and other relations to braided group coordinate
rings $BG_q$. The required relations correspond classically to the
characterisation of the group $G$ in $N\times N$ matrices. Then,
via the identification
\[ U(\cg_q)=B(R),\]
we make this same `geometrical' quotienting by braided
$q$-determinants etc. to get down from $U(\cg_q)$ to a version of
$U_q(\cg)$. Although it was known from \cite{FRT:lie} that
generators $l^+Sl^-$ in $U_q(g)$ obey the relations (\ref{BR}), the
prescription to make the required additional quotients when
starting from the quadratic algebra (\ref{BR}) is one of the
genuinely new results of braided group
theory\cite{Ma:exa}\cite{Ma:lie}. Also, the braided-Lie algebras
$\cg_q$ are again q-antisymmetric, although the general axioms for
such q-antisymmetry remain poorly understood. These remarks aim to
put the above anyonic Lie algebras into a wider context of
braided-Lie algebras and q-deformation.

In particular, the $su_2$ R-matrix gives $su_{2,q}$ as a
deformation of $u(1)\oplus su_2$, and $U(su_{2,q})=BM_q(2)$ is a
remarkable identification with  $2\times 2$-braided hermitian
matrices or $q$-Minkowski space in the general `twistor' or
spinorial-R-matrix approach of \cite{CWSSW:lor} (for the Lorentz
algebra) and \cite{Ma:exa}\cite{Ma:mec}\cite{Mey:lor} (for
spacetime algebras based on the relations (\ref{BR})). The
quotienting to $U_q(su_2)$ then corresponds geometrically to the
mass hyperboloid in $q$-Minkowski space.

\section{Anyfields}

There are probably many applications of anyonic Lie algebras and
anyspaces. One possible application is of course to consider the
anyonic variables as `organising' variables for anymultiplets. Thus
one can consider fields $\Phi(x,\theta)$ where $x$ is a usual
spacetime variable. The expansion
$\Phi(x,\theta)=\phi_0(x)+\theta\phi_1(x)+\cdots\theta^{n-1}\phi_{n-1}(x)$
gives the corresponding anymultiplet $(\phi_0,\cdots,\phi_{n-1})$.
Clearly a Lagrangian built from $\Phi$ could have hidden anyonic
properties not visible as a theory of fields $\phi_i$. Before
writing down such theories one probably needs the answer to the
following question:

\bigskip
{\em Problem 4} What is the anyonic analogue of the
super-Poincar\'e algebra?
\bigskip

Not knowing the answer to this does not, however, stop us from
proceeding to more `geometrical' anyfields (without knowing their
field theory). For example, anyonic gauge theory is introduced in
\cite{Ma:diag} as an example of braided group gauge theory. There
we consider $n=3$ and gauge fields $A(x,\theta)$. We refer to
\cite{Ma:diag} for details, but the gauge field corresponds to a
multiplet of 6 ordinary fields $(A_1,A_2,a_1,a_2,b_1,b_2)$ where
$A_1$ transforms as a $U(1)$ gauge field under the first component
gauge transformation (itself a multiplet $(c_1,c_2)$) and $A_2$ in a
more complicated way:
\eqn{anygauge}{A_1\mapsto A_1+\extd c_1,\quad A_2\mapsto A_2+\extd
 c_2+(1+e^{2\pi\imath\over 3})(A_1c_1-c_1A_1-c_1\extd c_1).}
The auxiliary fields $a_1,a_2,b_1,b_2$ also transform among themselves
and $A_1$. In fact, there are many variations of such gauge
theory according to the choices of differential calculi and gauge
group.

%\bibliographystyle{unsrt}
%\bibliography{biblio}

\end{document}